\documentclass[11pt]{article}
\usepackage{graphicx,amssymb,amsmath}
\usepackage{epsfig}
\usepackage{epsf}
\usepackage{times}
\usepackage{amsfonts}
\usepackage{latexsym}
\usepackage{bbm}     % blackboard characters (e.g. set of real numbers)
\usepackage{here}

\newtheorem{theorem}{Theorem}
\newtheorem{proposition}{Proposition}
\newtheorem{lemma}{Lemma}
\newtheorem{corollary}{Corollary}

\newtheorem{claim}{Claim}

\long\def\ignore#1{}

\newcommand{\later}[1]{{}}
\newcommand{\old}[1]{{}}

\newcommand{\NN}{\mathbb{N}} %  set of natural numbers
\newcommand{\RR}{\mathbb{R}} %  set of real numbers
\newcommand{\SH}{\mathbb{S}} %  set of unit vectors

\def\A{{\cal A}}
\def\B{{\cal B}}
\def\L{{\cal L}}

\def\U{{\cal U}}
\def\eps{\varepsilon}
\def\etal{{\it et~al.}\,}

\def\marrow{{\marginpar[\hfill$\longrightarrow$]{$\longleftarrow$}}}

\newcommand\comm[1]{\typeout{Used \string\comm...}\marrow\quad{\bf [[#1]]}\quad}

\title{On the number of tetrahedra with minimum,\\
  unit, and distinct volumes in three-space\thanks{A preliminary
    version of this paper appeared in the {\em Proceedings of the 18th
      ACM-SIAM Symposium on Discrete Algorithms (New Orleans, LA,
      2007)}, ACM Press, pp. 1114-1123.}}

\author{Adrian Dumitrescu\thanks{Department of Computer Science,
University of Wisconsin-Milwaukee,  WI 53201-0784, USA, email: {\tt
ad@cs.uwm.edu}}
\and Csaba D. T\'oth\thanks{Department of Mathematics,
MIT, Cambridge, MA 02139, USA, email: {\tt toth@math.mit.edu}}}

\begin{document}
\maketitle

\begin{abstract}
We formulate and give partial answers to several combinatorial
problems on volumes of simplices determined by $n$ points in
3-space, and in general in $d$ dimensions.
(i) The number of tetrahedra of minimum (nonzero) volume spanned by $n$
points in $\RR^3$ is at most $\frac{2}{3}n^3-O(n^2)$, and there are
point sets for which this number is $\frac{3}{16}n^3-O(n^2)$.
We also present an $O(n^3)$ time algorithm for reporting all
tetrahedra of minimum nonzero volume, and thereby extend an algorithm of
Edelsbrunner, O'Rourke, and Seidel. In general, for every $k,d\in
\NN$, $1\leq k \leq d$, the maximum number of $k$-dimensional
simplices of minimum (nonzero) volume spanned by $n$
points in $\RR^d$ is $\Theta(n^k)$.  (ii) The number of
unit-volume tetrahedra determined by $n$ points in $\RR^3$ is
$O(n^{7/2})$, and there are point sets for which this number is
$\Omega(n^3 \log \log{n})$.  (iii) For every $d\in \NN$, the minimum
number of distinct volumes of all full-dimensional simplices
determined by $n$ points in $\RR^d$, not all on a hyperplane, is
$\Theta(n)$.
\end{abstract}

\section{Introduction}

Typical Erd\H{o}s type problems in extremal discrete mathematics ask
for the minimum or maximum number of certain configurations over all
inputs of a given size. They are easy to formulate but often extremely
hard to answer. Their impact on mathematics and computer science has
been enormous, not only because of specific algorithms based on
combinatorial bounds but also because they have triggered the
development of theoretical and practical methods that turned out to be
applicable elsewhere.

Some of the most simply formulated yet notoriously hard Erd\H{o}s type
problems occur in combinatorial geometry. In 1946,
Erd\H{o}s~\cite{e-46} asked two questions on distances: (1) at most
how many times can a given distance occur among $n$ points in the
plane; (2) what is the minimum number of distinct distances determined
by $n$ points in the plane? The difficulty of the, so called, {\em
unit distance} and {\em distinct distance} problems is still to be
measured.
Erd\H{o}s and Purdy~\cite{ep-71,ep-76} generalized the unit- and
distinct distance problems to congruent (or repeated) simplices: What
is the maximum number of congruent $k$-dimensional simplices among $n$
points in $\RR^d$, for $1\leq k \leq d$?  No asymptotically tight
bound is known for this problem for $k\geq \lfloor d/2\rfloor$, (and
it is trivial for $1\leq k<\lfloor d/2\rfloor$ by Lenz-type
constructions).
Pach and Sharir~\cite{ps-92} derived bounds on the maximum number of
occurrences of the same angle
%(where each angle is represented by an {\em ordered} triple of points)
determined by a planar $n$-element point set; recently
Apfelbaum and Sharir~\cite{as-05} studied the analogous problem in
three-space.
A recent book by Bra\ss, Moser, and Pach~\cite{bmp-05} and a survey of
Pach and Sharir~\cite{ps-04} provide substantial details on these and
other similar problems.  In the sequel, we focus on problems about the
extremal number of simplices with certain volume properties in a point
set.

In 1967, A.~Oppenheim (see \cite{ep-95}) asked what is the maximum number of
unit area triangles determined by $n$ points in the plane.
Erd\H{o}s and Purdy gave an $O(n^{5/2})$ upper bound, and also
showed that a suitable section of the integer lattice yields
$\Omega(n^2 \log\log{n})$ such triangles ~\cite{ep-71}.
The currently best upper bound, $O(n^{44/19})$, due to
Dumitrescu, Sharir and Cs. T\'oth~\cite{dst-07},
recently improved an older $O(n^{7/3})$ bound of Pach and
Sharir~\cite{ps-92}.
Answering further questions of Erd\H{o}s and Purdy~\cite{ep-71},
Bra\ss, Rote, and Swanepoel~\cite{brs-01} showed the following two
results: (1) The maximum number of triangles of maximum area (or of
maximum perimeter) determined by $n$ points in the plane is exactly
$n$. (2) The maximum number of triangles of minimum (nonzero) area
determined by $n$ points in the plane is $\Theta(n^2)$.

In 1982, Erd\H os, Purdy, and Straus \cite{s-78,eps-82} considered the
generalization of the problem of distinct triangle areas to higher
dimensions and posed the following problem: Let $S$ be a set of $n$
points in $\RR^d$ not all in a hyperplane. What is the minimal number
$g_d(n)$ of distinct volumes of full-dimensional simplices with vertices
in $S$?
It is easy to see that $g_d(n) \leq \lfloor (n-1)/d \rfloor$ by
taking $d$ sets of about $n/d$ equally spaced points on parallel lines
through the vertices of an $(d-1)$-simplex. Erd\H{o}s, Purdy, and
Straus conjectured that equality holds at least for sufficiently large
$n$ (see also \cite{cfg-91}). Very recently, Pinchasi~\cite{p-07}
confirmed the conjecture in the plane after earlier work by Burton and
Purdy~\cite{bp-79} and Dumitrescu and Cs. T\'oth~\cite{dt-07}: The
minimum number of distinct (nonzero) triangle areas determined by $n$
noncollinear points is $\lfloor (n-1)/2\rfloor$.
In this paper we give a first partial result for
$d\geq 3$ by proving $g_d(n)=\Omega(n)$ for every $d\in \NN$, which is
optimal up to a multiplicative constant.
As expected, the three dimensional analogues of combinatorial
problems in the plane are often much harder. Only the recent few
years saw intensifying work on these problems in
three-space~\cite{as-05,aks-05,aps-04,fs-05,pps-04,sw-04}.

Algorithmic problems on sets of points in Euclidean space have been
often attacked using the duality transform and the machinery of
constructing hyperplane arrangements. The first applications of this
technique can be found in early works of Chazelle, Edelsbrunner,
Guibas, Lee, O'Rourke and Seidel~\cite{cgl-85,eg-89,ers-86}: among others,
the problem of computing a minimum area triangle in a given
set of $n$ points. A $O(n^2)$ time algorithm for such a task was given
independently in~\cite{cgl-85} and \cite{eg-89}. If there are
degeneracies in the set, the algorithm returns zero area and a triplet
of collinear points. An extension of the algorithm for finding a
minimum volume simplex among $n$ points in $\RR^d$ was given in
\cite{ers-86}: it runs in $O(n^d)$ time, but again reports zero volume
if there are $d+1$ points on a hyperplane.
Here we further extend the algorithms for two and three dimensions to
report all simplices of minimum {\em nonzero} volume, within the
same $O(n^2)$ and  $O(n^3)$ running times, respectively.

\paragraph{Our contribution.}
In this paper, we address classical problems on minimum, unit, and
distinct volume simplices.
We show that every set of $n$ points in three-space determines
$O(n^3)$ minimum (positive) volume tetrahedra, and there are points
sets in $\RR^3$ that span $\Omega(n^3)$ minimum volume tetrahedra.
Our techniques generalize to arbitrary dimensions:
For every $k,d\in \NN$, $1\leq k\leq d$, any set
of $n$ points in $\RR^d$ determines $O(n^k)$ minimum (positive) volume
$k$-dimensional simplices (each having $k+1$ vertices), and there are
points sets in $\RR^d$ that span $\Omega(n^k)$ minimum volume
$k$-dimensional simplices. In three-space, we also give a $O(n^3)$
time algorithm for reporting all tetrahedra of minimum nonzero volume,
and thereby extend an early algorithm of Edelsbrunner, O'Rourke, and
Seidel.
We prove that every set of $n$
points in $\RR^3$ determines at most $O(n^{7/2})$ unit-volume
tetrahedra, and there are point sets that span $\Omega(n^3 \log
\log{n})$ unit-volume tetrahedra.
Finally, we show that for every $d\in \NN$, any set of $n$ points in $\RR^d$,
not all in a hyperplane, determines at least $\Omega(n)$
full-dimensional simplices of distinct volumes, and there are point
sets for which this number is $O(n)$; this gives a first
answer to the question of Erd\H{o}s, Purdy, and Straus.

\paragraph{Organization.}
Section~\ref{S-toolbox} presents a collection of tools we use from
Euclidean geometry, previous results and extensions of previous results
adapted to our goals. We prove our main theorems on the number of
minimum-, unit-, and distinct volume tetrahedra in Sections
\ref{S-minimum}, \ref{S-unit}, and \ref{S-distinct}, respectively.
To simplify notation, we assume in our proofs that $\log{n}$ is an
integer.

\section{Toolbox} \label{S-toolbox}

We will frequently apply a classic result of Szemer\'edi and
Trotter~\cite{st-83} on the number of point-lines incidences in the
plane. (Since collinear points remain collinear under any affine
transformation, the result holds in arbitrary dimensions.) The
constant factor hidden in the asymptotic notation has been
significantly improved by the cutting method and by the theory of
crossing numbers~\cite{ps-04}; the current best constant is due to Pach et
al.~\cite{prt-04}. The Szemer\'edi-Trotter bound comes in two
equivalent formulations; we also state two immediate corollaries
that we use in our proofs. Given a point set $S$ in $\RR^d$, for any
integer $k\geq 2$, a line is called {\em $k$-rich} if it is incident
to at least $k$ points of $S$. We denote by $L_k$ the set of
$k$-rich lines.
\begin{theorem} {\rm (Szemer\'edi-Trotter~\cite{st-83})}.
\label{thm:st-83b}
Given $n$ points in $\RR^d$, $d\geq 2$, the
number of $k$-rich lines, $k\geq 2$, is
$$|L_k|=O\left(n^2/k^3+n/k\right). $$
\end{theorem}
\begin{theorem} {\rm (Szemer\'edi-Trotter~\cite{st-83})}.
\label{thm:st-83a}
The number of point-line incidences among $n$ points and $\ell$ lines
in $\RR^d$, $d\geq 2$, is
$$I(n,\ell)=O(n^{2/3}\ell^{2/3}+n+\ell). $$
\end{theorem}
\begin{corollary}\label{cor:st-83b}
  Given $n$ points in $\RR^d$, $d\geq 2$, the number of point-line
  incidences among the $n$ points and the $k$-rich lines, $k\geq 2$,
  is $I(n,|L_k|)=O(n^2/k^2+n)$.
\end{corollary}
Note that in Theorem~\ref{thm:st-83a}, the term $\ell$ is responsible
for the incidences on lines containing a single point.
\begin{corollary} \label{cor:st-83}
  Given $n$ points and $\ell$ lines in $\RR^d$, $d\geq 2$, each line
  containing at least two points, the number of point-line incidences
  among them is $ I(n,\ell)=O(n^{2/3}\ell^{2/3}+n)$.
\end{corollary}

Incidence bounds have been key components in many results in
combinatorial geometry. Some Szemer\'edi-Trotter-type bounds were
found for point-plane incidences in the space but they are either not
known to be tight or they hold in severely restricted settings only
(e.g., \cite{bk-03,ceg-90,et-05,ikr-06}). Instead of these
multidimensional bounds, we apply {\em planar} incidence bounds on the
{\em projections} of a finite point set $S\subset \RR^3$ onto planes in
certain directions determined by $S$.
The Szemer\'edi-Trotter bound on point-line incidences is tight in the
worst case, but does not hold for points with {\em multiplicities}.
Since the planar projection of a three-dimensional point set may have
an irregular distribution of multiplicities, we partition the point
sets (and the planar multiset) into subsets of roughly the same
multiplicities.
We demonstrate this technique in Section~\ref{S-unit}, where we
aggregate Szemer\'edi-Trotter bounds for various substructures and
give an upper bound on the number of incidences of lines and points
{\em with multiplicities} in the plane, and ultimately a bound on incidences
in three-space.

\subsection{Degenerate and nondegenerate planes}

We extensively apply the concepts of {\em degenerate} and {\em
  nondegenerate} planes in three-space. They were introduced recently
by Elekes and T\'oth~\cite{et-05}---here, we use the same concepts
with a different terminology. Consider a point set $S\subset \RR^3$.
For a constant $\alpha$, \ $0< \alpha \leq 1$, we say that a plane $A$
is {\em $\alpha$-degenerate} if at most $\alpha |S\cap A|$ points of
$S\cap A$ are collinear. In this paper, we fix $\alpha=\frac{2}{3}$
and say that a plane $A$ is {\em nondegenerate} if it is
$\frac{2}{3}$-degenerate; otherwise $A$ is {\em degenerate}. A
nondegenerate plane is always spanned by $S$, while a degenerate
plane may contain only one point, or collinear points. If $S$ lies in
a plane $A$, we refer to $S$ as degenerate or nondegenerate according
to the above definition.

We apply several consequences of a well known combinatorial geometric
result due to Beck~\cite{b-83}.

\begin{lemma} \label{lem:beck2} {\rm (Beck~\cite{b-83})}.
  Consider a set $S$ of $n$ points in the plane. If at most $\ell$
  points of $S$ are collinear, then $S$ determines at least
  $\Omega(n(n-\ell))$ distinct lines.
\end{lemma}

\begin{corollary} \label{cor:beck3}
  Consider a set $S$ of $n$ noncollinear points in the plane.  If at
  most $\ell$ points of $S$ are collinear, then $S$ determines at
  least $\Omega(n^2(n-\ell))$ nondegenerate triangles.
\end{corollary}
\begin{proof}
  If every line contains fewer than $n/2$ points, then by
  Lemma~\ref{lem:beck2}, $S$ determines $\Omega(n^2)$ line segments.
  If a line contains $n/2$ points, then there are $\Omega(n^2)$
  segments along this line. In either case, each of these line
  segments is the base for at least $n-\ell$ distinct nondegenerate triangles. Since
  we count every triangle at most three times, we obtain
  $\Omega(n^2(n-\ell))$ distinct triangles.
\end{proof}

For a set $S$ of $n$ points in the plane, let $u(S)$ denote the number
of unit area triangles determined by $S$. Let $\delta>0$ be a constant
such that $u(S)=O(n^{2+\delta})$ for every set $S$ of $n$ points in
the plane. It is conjectured that $\delta>0$ can be arbitrarily small;
we can assume $\delta\leq 1/3$ due to the bound $u(S)=O(n^{7/3})$ of
Pach and Sharir~\cite{ps-92} (the current best bound
$u(n)=O(n^{44/19})$ gives $\delta\leq 6/19$, cf. \cite{dst-07}).  We
use the general parameter $\delta$ when we further chisel the bound on
unit area triangles in the special case of degenerate planes.
We apply these results in Section~\ref{S-unit} in a charging scheme,
where we distribute unit area triangles spanned by a planar point set
among line segments in that plane.

\begin{lemma}\label{lem:unit}
  Let $S$ be a set of $n$ points in the plane, and $L$ be a line
  incident to exactly $\ell$ points of $S$. Then $S$ determines at
  most $O((n-\ell)\ell +(n-\ell)^{2+\delta})$ unit area triangles.
\end{lemma}
\begin{proof}
  Put $x=n-\ell=|S\setminus L|$, and let us denote the lines parallel
  to $L$ containing at least one point of $S\setminus L$ by
  $L_1,L_2,\ldots , L_q$, where $q$ is the number of such lines. For
  every $i=1,2,\ldots , q$, let $k_i=|S\cap L_i|$.
  We partition the unit area triangles into three subsets:
\begin{itemize}
\item[(i)]
$u_1$ counts triangles with all three vertices in $S\setminus L$;
\item[(ii)]
$u_2$ counts triangles with two vertices on $L$ and one in
$S\setminus L$;
\item[(iii)]
$u_3$ counts triangles with one vertex on $L$ and two in $S\setminus
L$;
\end{itemize}
(i) By the assumption, $u_1=t(S\setminus L)= O(x^{2+\delta})$.

\noindent (ii) For every line $L_i$, $i=1,2,\ldots, q$, there are at
most $k_i(\ell-1)$ unit area triangles $\Delta{abc}$, with $a,b \in L$
and $c\in L_i$. Therefore, $u_2\leq \sum_{i=1}^q k_i (\ell-1) \leq
\ell x$.

\noindent (iii) $S\setminus L$ determines less than $x^2/2$ line
segments $ab$. For segments $ab$ non-parallel to $L$, there are at
most two vertices $c\in L$ such that $\Delta{abc}$ has unit area
(because $c$ has to lie on one of two lines parallel to $ab$), so we
have less than $x^2$ unit area triangles of this kind.  For segments
$ab$ on a line $L_i$ parallel to $L$, there are at most $(k_i-1)\ell$
unit area triangles $\Delta{abc}$, with $a,b \in L_i$ and $c\in L$.
Summing over all lines $L_i$, there are at most $\ell \sum_{i=1}^q
(k_i-1) \leq \ell x$ triangles of unit area (since $\sum k_i \leq x$)
of this second kind.  Consequently $u_3=O(\ell x +x^2)$.
Altogether we have at most $\sum_{i=1}^3u_i= O(\ell x+x^{2+\delta})$
unit area triangles in $S$.
\end{proof}

\begin{corollary}\label{cor:dir} There is a constant $c_0>0$ with
  the following property. If a set $S$ of $n$ points in the plane
  determines $t$ nondegenerate triangles, then there is a subset
  $Q\subset S$ of $c_0t/n^2$ points and for each $p\in Q$ there is a
  set $\L_p$ of $c_0 n$ pairwise non-overlapping line segments spanned
  by $S$ such that
{\rm (i)} all segments in  $\L_p$ have a common endpoint $p$
% one endpoint of each segment in  $\L_p$ is $p$
and {\rm (ii)} the sets of segments $\L_p$, $p\in Q$, are pairwise disjoint.
\end{corollary}

\begin{proof}
  First assume that $S$ is nondegenerate. By Corollary \ref{cor:beck3},
  we have $t=\Omega(n^3)$ and $S$ spans $\Omega(n^2)$ distinct lines.
  For any $k$, $1\leq k\leq \sqrt{n}$, at most $O(n^2/k^3)$ lines are incident to more than
  $k$ points by Theorem~\ref{thm:st-83b}. There is a constant $c_1>0$
  such that the set $\L_1$ of lines incident to at least 2 and at most
  $c_1$ points still contains $\Omega(n^2)$ lines. For each line
  $L\in \L_1$, choose an arbitrary point $p$ among the at most $c_1$ points
  in $S\cap L$ and place $L$ into the set $\L_p$. The resulting sets $\L_p$ are
  disjoint and each contain at most $n-1$ lines. For at least
  $\Omega(n)$ points $p\in S$, the set $\L_p$ contains at least
  $\Omega(n)$ lines.

  Next assume that $S$ is degenerate, that is, there is a line $L$
  incident to $\ell \geq \frac{2}{3}n$ points. In this case, every
  point off $L$ forms $\Omega(n^2)$ triangles with point pairs from
  $L$, and so $t=\Omega(n^2(n-\ell))$. Each of the $n-\ell$ points off
  the line $L$ is incident to $\Omega(n)$ lines spanned by $S$. For
  any constant $k$, $0<k\leq \sqrt{n-\ell}$,
  at most $O((n-\ell)^2/k^3)$ lines are incident to more than $k$ points
  in $S\setminus L$ by Theorem~\ref{thm:st-83b}. There is a constant
  $c_2>0$ such that the set $\L_2$ of lines incident to at least 2 and
  at most $c_2$ points contains $\Omega(n(n-\ell))$ elements. For each line
  in $\L_2$, choose an arbitrary incident point $p$ and assign the line to
  $\L_p$. For at least $\Omega(n-\ell)$ points $p\in S\setminus L$,
  the resulting set $\L_p$ contains at least $\Omega(n)$ lines.
\end{proof}

\begin{corollary}\label{cor:unit+}
  If a set of $n$ points in the plane determines $t$ nondegenerate
  triangles, then at most $O(tn^{\delta-1})$ of them have unit area.
\end{corollary}

\begin{proof}
  First, suppose that $S$ is nondegenerate. By Lemma~\ref{lem:beck2},
  $t=\Omega(n^3)$, hence $u(S)=O(n^{2+\delta})=O(tn^{\delta -1})$.

  Next, suppose that $S$ is degenerate. A line $L$ passes through
  $\ell \geq 2n/3$ points. Put $x=|S\setminus L|=n-\ell$. As in the
  proof of Corollary \ref{cor:dir}, we have $t=\Omega(n^2 x)$. It
  follows, using Lemma~\ref{lem:unit}, that $u(S)=O(x
  \ell+x^{2+\delta}) = O(xn^{1+\delta})=O(tn^{\delta -1})$.
\end{proof}

\subsection{Beck's lemma in higher dimensions}

Beck~\cite{b-83} also showed that Lemma~\ref{lem:beck2} can be extended
to higher dimensions:

\begin{lemma} {\rm (Beck~\cite{b-83})}. \label{lem:b}
For any $d\geq 2$ there exist constants $\beta_d,\gamma_d\in(0,1/2]$
such that, for any set $S$ of $n$ points in $\RR^d$, at least one of
the following holds.
\begin{enumerate}
\item a hyperplane contains more than $\beta_d n$ points of $S$; or
\item the $d$-tuples of $S$ span at least $\gamma_d n^d$ distinct
hyperplanes.
\end{enumerate}
\end{lemma}

\section{Minimum volume tetrahedra} \label{S-minimum}

We show that $n$ points determine at most $2n^3/3-O(n^2)$ minimum
volume tetrahedra in three space. The upper bound is based on a new
charging scheme which assigns each tetrahedron of minimum volume to
one of its four faces. We then extend our charging scheme
and show that for any fixed $k,d\in \NN$, $1\leq k\leq d$, the number of
$k$-dimensional simplices of minimum (nonzero) volume in $\RR^d$ is
$O(n^k)$, where the constant of proportionality depends only on $k$
and $d$. This bound is best possible apart from the constant factor.
\begin{theorem} \label{T-min}
  The number of tetrahedra of minimum (nonzero) volume spanned by $n$
  points in $\RR^3$ is at most $\frac{2}{3}n^3-O(n^2)$, and there are
  point sets for which this number is $\frac{3}{16}n^3-O(n^2)$. Given
  $n$ points in $\RR^3$, all tetrahedra of minimum nonzero volume can
  be reported in $O(n^3)$ time and $O(n^2)$ working space.
\end{theorem}
\paragraph{Lower bound.}
The lower bound construction is simple.
Form a rhombus with unit sides in the $xy$-plane from two equilateral
triangles $\Delta{abc}$ and $\Delta{abe}$ with a common side $ab$. 
Extend it to a prism in 3-space, and  make a unit volume tetrahedron with
vertices $a$, $b$, $c$, and $d$, on the four vertical lines. Replace each
of $a$, $b$, $c$, and $d$, with $n/4$ equally spaced points with inter-point
distances $\eps$ along these lines, for a sufficiently small $\eps>0$
(assume $n$ is divisible by 4). Observe that each tetrahedron with one
vertex on each of the four lines has volume very close to 1,
but the minimum volume is $O(\eps)$, given by tetrahedra with two
consecutive vertices on a line, and two vertices on any two of the
other three lines. The number of such tetrahedra is
${4\choose 3}3(n/4-1)(n/4)^2 =\frac{3}{16}n^3-O(n^2)$.
All other tetrahedra have zero volume.

\paragraph{Upper bound.}
Let $S$ be a set of $n$ points in
$\RR^3$. Denote by $T=T(abcd)$ the tetrahedron determined by four
non-coplanar points $a,b,c,d \in S$ and by ${\rm vol}(T)$ its volume.
Similarly, let $T(s_1,s_2)$ denote the tetrahedron determined by the
endpoints of two line segments $s_1$ and $s_2$. The key to our upper
bound is the following charging scheme: assign every (nondegenerate)
tetrahedron determined by $S$ to one of its four faces as follows.
Assign $T(abcd)$ to a triangle face of maximum area among the
faces adjacent to a diameter of $T(abcd)$. We show that at most a
constant number of minimum volume tetrahedra are charged to every
triangle, and this yields the desired $O(n^3)$ upper bound.

\begin{figure}[htbp]
\centerline{\epsfxsize=3.9in \epsffile{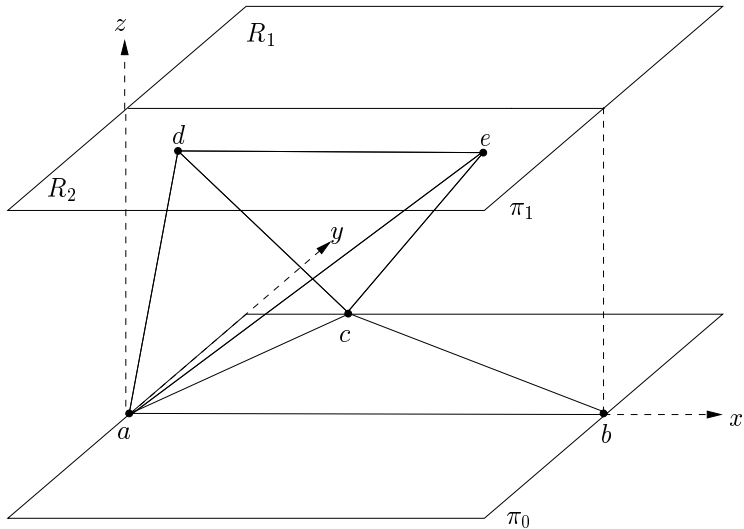}}
\caption{\small Illustration to the proof of Theorem~\protect{\ref{T-min}}.
Planes $\pi_0$ and $\pi_1$, and an example for
the second case in the argument ($de$ is parallel to $ab$).}
\label{fig:min}
\end{figure}

Consider a (nondegenerate) triangle $\Delta abc$ with $a,b,c\in S$.
Let $ab$ be a diameter of $\Delta abc$. Choose a coordinate system
such that $a$ is the origin, $ab$ lies on the $x$-axis, and $\Delta
abc$ lies in the $xy$-plane. Refer to Fig.~\ref{fig:min}.

Let $|ab|=x_0$, and denote by $y_0$ the height from $c$ in $\Delta
abc$. Assume that the minimum volume of a tetrahedra spanned by $S$
is $v_0$. All points $d$ with ${\rm vol}(T(abcd))=v_0$ must lie in
two horizontal planes at distance $z_0=3v_0/{\rm area}(\Delta abc)$
from  the $xy$-plane. Let us consider the plane $\pi_1:z=z_0$ for
now. By our assignment, all points $d$ must lie in the interior of
an axis-aligned rectangle $R=(0,x_0)\times (-y_0,y_0)\times [z_0]$
in the plane $\pi_1$ (otherwise $ab$ would not be the diameter or
$\Delta abd$ would have larger area than $\Delta abc$; note that
point $d$ cannot lie on the boundary of $R$). Partition the open
rectangle $R$ into two rectangles $R_1=(0,x_0)\times (-y_0,0]\times
[z_0]$ and $R_2=(0,x_0)\times (0,y_0)\times [z_0].$

\begin{claim}
For $i=1,2$, there is at most one point $d\in S\cap R_i$ such that
${\rm vol}(T(abcd))$ is minimum.
\end{claim}
\begin{proof}
Assume there are two points $d,e\in S\cap R_i$  for some $i \in
\{1,2\}$. We will pick a side $s\in \{ab,ac\}$ so that $s$ is not parallel
with $de$. Then $T(s,de)$ is a nondegenerate tetrahedron. We show
that ${\rm vol}(T(s,de)) < v_0$.

Assume first that $de$ is not parallel to $ab$, and pick $s=ab$. The
segment $de$ lies in rectangle $R_i$, $i \in \{1,2\}$; observe that
at most one element of $\{d,e\}$ may lie on the boundary of $R_i$.
Denote by $\ell_1$ and $\ell_2$ the two lines containing the two
sides of $R_i$ parallel to the $x$-axis. The volume of $T(s,de)$
strictly increases if we move $d$ to position $d'$ and $e$ to
position $e'$ along the line $de$ such that $d'\in \ell_1$ and
$e'\in \ell_2$. The volume of $T(ab,d'e')$ is $x_0y_0z_0/6=v_0$,
hence ${\rm vol}(T(ab,de)) < v_0$ (e.g., the base triangle of
$T(ab,d'e')$ lies in the $xz$-plane, and its height is $y_0$).
A contradiction.

Next assume that $de$ is parallel to $ab$, and pick $s=ac$. In this
case, $d$ and $e$ are strictly in the interior of $R_i$. If we
replace $d$ and $e$ by points $d'$ and $e'$ on the line $de$ such
that $d'e'$ has the same $x$-extent as the rectangle $R_i$, then
${\rm vol}(T(ac,de)) < {\rm vol}(T(ac,d'e'))$. Draw two lines
parallel to $s=ac$ through $d'$ and $e'$. Let these two
lines intersect the $xz$-plane at points $d''$ and $e''$. (Here we
use the fact that if the vertices of a tetrahedron are on three
parallel lines, one can shift the single points along the
corresponding lines, and the volume remains the same.) Since we
moved $d'$ and $e'$ along lines parallel to $s$, the volume of
$T(s,d'e')$ is the same as that of $T(s,d''e'')$, and this volume is
$x_0y_0z_0/6=v_0$ (e.g., the base triangle of $T(ac,d''e'')$ lies
in the $xz$-plane, and its height is $y_0$). A contradiction.
\end{proof}

Symmetrically, the plane $z=-z_0$ also contains at most 2 points $d$
with ${\rm vol}(T(abcd))=v_0$, hence the number of minimum volume
tetrahedra is at most $4 {n \choose 3} =\frac{2}{3}n^3 -O(n^2)$.

\paragraph{Reporting minimum nonzero volume tetrahedra.}
Next, we devise an algorithm for counting and reporting all minimum
(nonzero) volume tetrahedra spanned by a given set $S$ of $n$ points
in three-space.

A plane $h$ and a point $p\not\in h$ form a {\em slab}, which is
defined as the open region bounded by $h$ and the plane $h'$, which
is parallel to $h$ and incident to $p$. We say that the slab is {\em
empty}, if it is disjoint from $S$. Bra\ss, Rote, and
Swanepoel~\cite{brs-01} observed that if $\Delta abc$ is a minimum
area triangle in the plane, then the relative interior of the line
segment $ab$ and the slab formed by the line through $ab$ and the
point $c$ is empty of $S$. This observation readily generalizes to
three-space: If $abcd$ is a minimum volume tetrahedron, then the
triangle $\Delta abc$ is a minimum area triangle in the plane
$A(abc)$ spanned by $\Delta abc$, and the open slab formed by
$A(abc)$ and point $d$ is empty.

Let $\A$ denote the set of planes spanned by $S$. We compute for
every empty slab formed by a plane $h\in \A$ and point $p$, the
number of tetrahedra whose base is a minimum area triangle in $h$
and whose 4th vertex lies in a plane $h'$ parallel to $h$ and
incident to $p$; we also record the volume of these tetrahedra
(which might not be the global minimum). We then sum up these
quantities for those empty slabs where this volume is {\em the
minimum}. We can compute the number of minimum volume tetrahedra,
since each of them is counted exactly four times, once for each face.
If we record all minimum area triangles in every plane $h$ and every point
in plane $h'$, then we can also {\em report} all
minimum volume tetrahedra (we have already shown that there are only
$O(n^3)$ of them).

Our algorithm works on the dual arrangement $S^*$.
Consider the (well-known) duality transform which maps a point
$p(a,b,c)\in S$ to the (nonvertical) plane $p^*$ with equation
$p^*:z=ax+by-c$. Conversely, a nonvertical plane $q$ with equation
$z=ax+by+c$ is mapped to the point $q^*=(a,b,-c)$.
Duality $D$ satisfies $(p^*)^*=p$, for any point $p$, and
$(q^*)^*=q$, for any nonvertical plane $q$. The duality preserves
point-plane incidences and it reverses the above-below relationship,
which is understood with respect the the $z$-axis. Furthermore, it
preserves vertical distances for point-plane pairs. More precisely,
point $p$ is above plane $q$, if and only if point $q^*$ is above
plane $p^*$, and the vertical distance between the point and the
plane in each pair is the same.
Every plane $h$ spanned by the point set $S$ corresponds to a vertex
$h^*$ of the arrangement. An empty slab determined by $h$ and a
point $p$ corresponds to a vertical line segment between $h^*$ and
plane $p^*$ that does not pierce any plane of the arrangement $S^*$.

First we solve a planar problem: count the minimum (nonzero)
area triangles in a plane efficiently. For a plane $h\in \A$, let
$n_h=|S\cap h|$, let $\ell_h$ be the number of lines spanned by
$S\cap h$, and let $M_h$ denote the number of minimum (nonzero) area
triangles determined by $h\cap S$. Let $\L$ denote the set of lines
spanned by $S$.
\begin{proposition}\label{pro:time}
Consider a plane $h\in \A$. Given the list of $n_h=|S\cap h|$ points
in $h$, and for every line $L\in \L$ the length and the number of
shortest (that is, {\em minimum length}) segments in $S\cap L$,
we can compute $M_h$ in $O(n_h \ell_h)$ time.
\end{proposition}
\begin{proof}
For every line $L\in \L$, $L\subset h$, compute the number of minimum area
triangles in the empty slabs formed by $L$ (that is, triangles such
that one side is a shortest segment $S\cap L$ and a third vertex
is a closest point to $L$). Since the number of shortest
segments along $L$ is given, it suffices to check the distance of
all points of $S\cap h$ to the line $L$. We can enumerate these
triangles in $O(n_h)$ time. Summing this over all $\ell_h$ lines,
$M_h$ can be computed in $O(n_h \ell_h)$ time.
\end{proof}

\noindent We use the Szemer\'edi-Trotter Theorem (Theorem~\ref{thm:st-83b}) to
bound the time to compute $M_h$ for all planes $h\in \A$.

\begin{proposition}\label{cor:n3}
With the notation defined above, we have $\sum_{h\in \A} n_h \ell_h
= O(n^3)$.
\end{proposition}
\begin{proof}
Denote by $\L_i$ the set of lines incident to at least $2^i$ but
fewer than $2^{i+1}$ points of $S$, for $i=1,2,\ldots , \log n$. By
Corollary~\ref{cor:st-83},
we have $\sum_{h\in \A} n_h \ell_h
\leq \sum_{L\in \L} |L\cap S| \cdot |S\setminus L|
\leq n \sum_{L\in \L} |L\cap S|
%\leq n \sum_{i=1}^{\log n} |\L_i| 2^{i+1}
\leq n \sum_{i=1}^{\log n} (n^2/2^{2i}+n)
= O(n^3+ n^2\log n)=O(n^3).$
\end{proof}

We can now count the minimum (nonzero) volume tetrahedra spanned by $S$.
Assume that the entire dual arrangement of $S^*$ and its incidence
graph is available; it can be computed in $O(n^3)$ time, and stored in
$O(n^3)$ space~\cite{as-00,e-87}. The arrangement is a cell complex
with faces of dimension 0, 1, 2, and 3.  For every vertex, we store
the incident planes, and the incident ridges (lines formed by the
intersection of two planes); for every ridge, we store the incident
vertices (sorted along the ridge) and planes of the arrangement
(sorted around the ridge).

We perform the following three tasks: (1). Preprocess the lines
spanned by $S$. For every line $L\in \L$, compute the number of
shortest segments in $S\cap L$ and record their (common) length.
The sorted list of points along each line can be extracted from the
arrangement. Since there are at most ${n\choose 2}$ lines, each
containing no more than $n$ points, this can be done in $O(n^3)$ time.
(2) For every vertex $h^*$, compute the number $M_h$ of minimum
(nonzero) area triangles lying in the plane $h$. By Proposition~\ref{pro:time},
this can be done in $O(n_h\ell_h)$ time for a plane $h$. By
Proposition~\ref{cor:n3}, the total time over all planes in $\A$ is
$\sum_{h\in \A} O(n_h\ell_h)= O(n^3)$. The list of points in each
plane $h$ is also available from the arrangement. (3) Sweep the dual
arrangement $S^*$ with a horizontal plane. For every vertex $h^*$,
drop a vertical line through $h^*$ and find the face $f$ of the cell hit
by this line, and the intersection point $x$ of the vertical line
with $f$: the face $f$ can be of dimension 0, 1, or 2.
Again, the number $N_f$ of planes in $S^*$ incident to $f$ can be
extracted from the arrangement. The empty slab formed by $h$ and
$x^*$ determines $M_h N_f$ tetrahedra with a minimum area base
triangle in $h$ and a 4th point in plane $x^*$. For every vertex
$h^*$ of the arrangement, we record $M_h N_f$ and the volume of the
corresponding tetrahedra. If this volume turns out to be the minimum
volume, all these $M_h N_f$ tetrahedra can also be reported in a
second pass over the data. Since the total number of minimum volume
tetrahedra is $O(n^3)$, all reporting takes $O(n^3)$ time.

If the entire dual arrangement $S^*$ and its incidence graph is
available, then our algorithm requires $O(n^3\log n)$ time and
$O(n^3)$ space in the real RAM model of computation (where algebraic
operations over reals have unit cost). The space requirement can be
reduced to $O(n^2)$ if we use a plane sweep algorithm to perform
tasks (2) and (3).
It sweeps the arrangement with a horizontal plane and
keeps in memory the cells intersecting the sweep plane, whose total
complexity is $O(n^2)$ by the zone theorem~\cite{ess-93}. When the
sweep plane arrives at a vertex $h^*$, it can provide the set of
planes and ridges incident to $h^*$ and so we can perform task (2).
The plane sweep algorithm takes $O(n^3 \log n)$ time (since it
maintains an event queue, which is updated in $O(\log n)$ time for
each vertex of the arrangement). The time can be reduced to $O(n^3)$
while maintaining an $O(n^2)$ working space by the {\em topological sweep}
method of Edelsbrunner and Guibas~\cite{eg-89} (originally designed
for line arrangements in the plane, and adapted to three-dimensions by
Anagnostou {\em et~al.}~\cite{agp-90}).

Edelsbrunner and Guibas~\cite{eg-89} and Chazelle
\etal~\cite{cgl-85} gave quadratic time sweep-line algorithms for
the problem of finding a triangle of minimum (possibly zero) area.
It is easy to rewrite our 3D algorithm for the planar problem and
report all triangles of minimum nonzero area determined by a set of
$n$ points in $O(n^2)$ time and space. When the sweep line arrives at vertex
$v$ of the dual arrangement, pairs of lines incident to $v$
 with consecutive slopes correspond to minimum length segments
determined by the points of $S$ whose duals intersect at $v$. By
drawing vertical segments at each vertex of the arrangement, the
algorithm examines all triangles whose base is a minimum length
segment on some line of $\L$, and the third vertex is closest to the
line through the base.

\paragraph{Generalization to $k$-simplices in $d$-space}

Theorem~\ref{T-min} can be generalized to arbitrary dimensions.

\begin{theorem} \label{D-min}
  For every (fixed) $k,d\in \NN$ with $1\leq k\leq d$, the number of
  $k$-dimensional simplices of minimum (nonzero) volume spanned by
  $n$ points in $\RR^d$ is $O(n^k)$, and there are point
  sets for which this number is $\Omega(n^k)$.
\end{theorem}
\begin{proof}
{\bf Lower bound.}
Consider $k$ parallel lines such that any two are at unit distance
apart from each other; place $n/k$ equally spaced
points on each line (assume $n$ is a multiple of $k$). Every
$k$-simplex with positive volume has two
vertices on one line and one vertex on each of the other lines. If the
$k$-simplex has minimum volume, the vertices along the same line must
be consecutive.  Conversely, every $k$-simplex with two consecutive
vertices on one line and one arbitrary point on each of the other two
lines has the same (minimum) volume. The number of such $k$-simplices is
$k(n/k-1)(n/k)^{k-1} = k^{1-k}n^k-O(n^{k-1})=\Omega(n^k)$.

\paragraph{Upper bound.}
Let $S$ be a set of $n$ points in $\RR^d$. Assign every
(nondegenerate) $k$-simplex $T$ determined by $S$ to one of its $k+1$
faces as follows. We choose the $k$ vertices of a face of $T$ in $k$
steps. First choose an arbitrary vertex $a_1$ of $T$. In step
$i=2,3,\ldots ,k$, choose a vertex $a_i$ of $T$ which is at maximum
(Euclidean) distance from the $(i-2)$-dimensional affine subspace
$A_{i-1}$ spanned by $\{a_1,a_2,\ldots , a_{i-1}\}$.
We show that at most a constant number of minimum volume
$k$-simplices are assigned to every $(k-1)$-simplex spanned by $S$.
This yields an $O(n^k)$ upper bound.

Assume that the minimum volume of a $k$-simplex is 1.
Consider a (nondegenerate) $(k-1)$-simplex $F=\{a_1,a_2,\ldots , a_k\}$
such that $a_i$, $i=2,3,\ldots,k$, is a furthest point in $\{a_i,a_{i+1},\ldots,a_k\}$
from the affine subspace $A_{i-1}$.
Choose an orthogonal coordinate system such that $a_1$ is the origin and $A_i$,
$i=2,3,\ldots ,k$, is the $(i-1)$-dimensional subspace spanned by the
first $i-1$ coordinate axes. Denote the $k-1$ extents of the axis-aligned
bounding box of $F$ by $y_1,y_2,\ldots y_{k-1}$; and choose $y_k$ so
that $\Pi_{i=1}^k y_i = k!$.
So $k!$ is the volume of the bounding box of any unit volume $k$-simplex
with a $(k-1)$-dimensional face $F$ in the coordinate system described.
By the choice of the points $a_1,a_2,\ldots ,a_k$, we have
$y_1\geq y_2\geq \ldots \geq y_{k-1}$. If our charging scheme assigns
a minimum volume $k$-simplex $F\cup \{b\}$ to $F$, then $y_{k-1} \geq y_k$,
since we chose $a_k$ such that it is at least as far from $A_{k-1}$ as $b$.

Every point $b$ such that $F\cup \{b\}$ is a $k$-simplex of a given
volume lies in a hypersurface $M$, which is the Minkowski sum $M=A_k+y_k
\SH^{d-k}$ of the $(k-1)$-dimensional Euclidean space $A_k$ spanned by
$F$ (in the first $k-1$ coordinates) and a sphere $y_k\SH^{d-k}$ of radius
$y_k$ (in the remaining $d-k+1$ coordinates). Every point $b\in M$ has a
unique decomposition into the vector sum $b=b_f+b_s$ with $b_f\in A_k$
and $b_s\in y_k\SH^{d-k}$.

If the $k$-simplex $F\cup \{b\}$, $b\in S \cap M$, is assigned to $F$, then
$b$ must satisfy $k-1$ constraints: $b$ is at most as far from the
subspace $A_{i-1}$ as $a_i$ for $i=2,3,\ldots , k$. These constraints
are satisfied if the $i$-th coordinate of $b_f$ is at most $y_i$,
$i=1,2,\ldots , k-1$. That is, $b_f$ lies in an axis-aligned rectangle
$R\subset \RR^{k-1}$ of extents $2y_i$, for $i=1,\ldots,k-1$.

Partition $R$ into $\lceil \sqrt{8d}\rceil^{k-1}$ congruent
axis-aligned rectangles of extents $2y_i/\lceil \sqrt{8d}\rceil$, each
for $i=1,2,\ldots,k-1$. Partition $\SH^{d-k}$ into $O(d^{(d-k)/2})$
convex regions of diameter $y_k/\sqrt{2d}$ each.
These two partitions give a partition of the Minkowski sum
$R+y_k\SH^{d-k}$ into
$O((8d)^{(k-1)/2} \cdot d^{(d-k)/2}) \leq O(2^{3d/2} \cdot d^{d/2}) \leq
O(d^{2d})$ regions. Next we show that each region
contains at most one point of $S$. This confirms that at most
$O(d^{2d}) = O(1)$ minimum $k$-simplices are assigned to $F$, completing our
proof.

Suppose, by contradiction, that $b,c\in S$ are assigned to $F$ and both lie
in the same region of $R+y_k\SH^{d-k}$. Since $b,c\in M$, both $F\cup
\{b\}$ and $F\cup \{c\}$ are $k$-simplices of minimum volume. We show
that $F\cup \{b,c\}$ spans a $k$-simplex of a strictly smaller positive
volume, which is a contradiction. Consider the maximum
index $j\in \{2,\ldots , k\}$ for which the set $\{a_1,\ldots ,
a_{j-1}, b,c\}$ is a
nondegenerate $j$-simplex. (Such an index $j\in \{2,\ldots , k\}$ exits,
since $\{a_1,b,c\}$ is a noncollinear triple.) We show that the $k$-simplex
$F\cup \{b,c\}\setminus \{a_j\}$ has smaller volume than $F\cup \{b\}$.
Note that it is enough to compare the volumes of  the $j$-simplices
$T_1=\{a_1,\ldots,a_j,b\}$ and $T_2=\{a_1,\ldots,a_{j-1},b,c\}$.
To compare the volumes of $T_1$ and $T_2$, it is enough to compare the
distances of $a_j$ and $c$ to the affine subspace $A'$ spanned by
$\{a_1,a_2,\ldots, a_{j-1},b\}$, that is, ${\rm dist}(a_j,A')$ and
${\rm dist}(c,A')$.
If we shift the line segment $bc$ in an affine
subspace parallel to $A_j$, then the volumes of  $T_1$ and $T_2$ do not change,
so we may assume that segment $a_1b$ is orthogonal to $A_j$.
Then, we have ${\rm dist}(a_j,A') = y_j\geq y_k$.
The first $j-1$ coordinates of the vector 
$\overrightarrow{bc}$ are parallel to $A_{j-1}\subset A'$, 
and so they can be ignored when computing ${\rm dist}(c,A')$.
Each of the coordinates $j,j+1,\ldots , k-1$ of $\overrightarrow{bc}$
is at most a $1/\sqrt{2d}$ fraction of the corresponding coordinate of $a_j$; 
the $k$-th coordinate of $\overrightarrow{bc} \subset M$ is $0$, 
and each of the last $d-k$ coordinates of $\overrightarrow{bc}$
is at most $y_k/\sqrt{2d}$. Hence
$${\rm dist}(c,A')\leq \sqrt{\sum_{i=j}^{k-1}
\left(\frac{y_i}{\sqrt{2d}}\right)^2 + \sum_{i=k+1}^d
\left(\frac{y_k}{\sqrt{2d}}\right)^2}
<\frac{y_j}{\sqrt{2}}
<y_j = {\rm dist}(a_j,A').$$
This confirms that $c$ is strictly closer to $A'$ than $a_j$, hence
$0<{\rm vol}(T_2)<{\rm vol}(T_1)$.
\end{proof}

\section{Unit volume tetrahedra} \label{S-unit}

\paragraph{Lower bound.}
As mentioned in the introduction, Erd\H{o}s and Purdy showed that a
suitable section of the integer lattice has $\Omega(n^2 \log\log{n})$
unit-area triangles \cite{ep-71}. This immediately gives a lower bound
of $\Omega(n^3 \log\log{n})$ unit volume tetrahedra (e.g., by placing
two such lattice sections, with $n/2$ points each, in two parallel planes).

\paragraph{Upper bound.}
Pach and Sharir~\cite{ps-92} showed
that every point $a$ of a set $S$ of $n$ points in the plane is
incident to at most $O(n^{4/3})$ unit area triangles spanned by $S$.
This bound follows easily from the Szemer\'edi-Trotter bound on the
number point-line incidences: For every point $b\in S\setminus
\{a\}$, there is a unique line $L_b$ parallel to $ab$ such that
every point $c\in S\cap L_b$ above $L_b$ gives a unit area triangle
$\Delta abc$; there are $O(n^{4/3})$ point line incidences between
$S$ and the lines $L_b$. This bound is tight, since the
Szemer\'edi-Trotter bound is tight: there are point sets with
$\Omega(n^{4/3})$ unit triangles having a common vertex: Place a set
$S_1$ of $n/2$ points and a set $\L_2$ of $n/2$ lines with
$\Omega(n^{4/3})$ incidences in the plane; $\L_2$ determines a point
set $S_2$ with $\L_2=\{L_b:b\in S_2\}$, and so $S=S_1\cup S_2\cup
\{a\}$ spans $\Omega(n^{4/3})$ unit area triangles incident to $a$.

We follow a similar strategy in three dimensions. Instead of
bounding directly the number of unit-volume tetrahedra with a common
vertex or a common edge, we design a charging scheme. We assign
every unit volume tetrahedron to a line segment lying in the plane
containing one of its four faces. Ideally, a tetrahedron is assigned
to one of its edges, but this does not always hold in our charging
scheme. We show that at most $O(n^{3/2})$ tetrahedra are assigned to
every segment; which immediately gives an $O(n^{7/2})$ bound on the
number of unit volume tetrahedra.

\begin{theorem}\label{thm:unit+}
The number of unit-volume tetrahedra determined by $n$ points in
$\RR^3$ is $O(n^{7/2})$.
\end{theorem}

\begin{proof}
Let $S$ be a set of $n$ points in $\RR^3$, and $\A$ denote the set
of all planes spanned by $S$. We follow the convention that in a
tetrahedron $t=abcd$, vertex $d$ lies above the plane $A(abc)$
spanned by $\Delta abc$. For every triangle $\Delta abc$, every
vertex $d$ for which $abcd$ has unit volume lies in a plane $B(abc)$
parallel to plane $A(abc)$. We decompose the set $\U$ of unit volume
tetrahedra into two subsets:
\begin{itemize}\topsep=0pt\itemsep=-3pt\parsep=1pt
\item
$\U_1$ contains every $t\in \U$ where $B(abc)$ is nondegenerate;
\item
$\U_2$ contains every $t\in \U$ where $B(abc)$ is degenerate.
\end{itemize}
{\bf Consider $\U_1$.} For two distinct points $v,w\in S$, let
$\A_{vw}$ denote the set of planes spanned by $S$ containing the
segment $vw$.

\paragraph{Assigning triangles to segments.}
We design a charging scheme where we assign every triangle $\Delta
abc$ to a segment in the plane $A(abc)$. This induces a charging
scheme for unit volume tetrahedra: If $\Delta abc$ is assigned to
segment $vw$, then we assign every unit volume tetrahedron $abcd$ to
segment $vw$, too.

The assignment is done in each plane $A\in \A$ independently.
Consider a plane $A\in \A$. We proceed in two stages: in the first
stage we compute, for every triangle $\Delta abc$, a
collection $E(abc)$ of segments in $A(abc)$; in the second stage, we
assign $\Delta abc$ to one of the segments in $E(abc)$. The assignment
is fairly elaborate; intuitively, we pursue two goals:
(1) every segment $vw$ should be assigned to ``few'' triangles of the same
area, and (2) a  segment $vw\subset A(abc)$ should {\em not} be assigned to
a triangle $\Delta abc$ if the plane $B(abc)$ contains a
``rich'' line parallel to $vw$. We continue with the details.

Denote by $k_A= |S\cap A|$ the number of points in $A$, and by $t_A$
the number of triangles spanned by $S\cap A$.
Corollary~\ref{cor:dir} tells us that there is a set $Q_A\subset
S\cap A$ of $c_0 t_A/k_A^2$ points, for some absolute constant
$c_0>0$, such that each $p\in Q_A$ is incident to a set $\L_p$ of $c_0
k_A$ line segments determined by $S\cap A$, and the line segments in
$\L_p$, $p\in Q_A$, are all distinct. Denote the collection of all
these segments by $E_A = \bigcup_{p\in Q_A}\L_p$.
We have $|E_A|=c_0^2 t_A/k_A$.
For every triangle $\Delta abc\subset A$, let $F(abc)$ denote the set of
segments in $A$ parallel to any one of the $c_0 k_A/2$ richest lines
in the plane $B(abc)$. Let $E(abc):= E_{A(abc)}\setminus F(abc)$ be the set
of segments to which $\Delta abc$ may be assigned. Since $|\L_p
\setminus F(abc)|\geq c_0 k_A/2$ for every $p\in Q_A$, we have
$|E(abc)|\geq |E_A|/2 = c_0^2 t_A/(2k_A)$.

For two points $v,w\in S\cap A$, let $k_{vw}(A)$ denote the number
of points in $S\cap A$ that do not lie on the line through $vw$. If
$vw\in E_A$, then one endpoint of $vw$ is in $Q_A$, and by
Corollary~\ref{cor:dir}, we have $k_{vw}(A) = \Theta(k_A)$. Hence,
we have $|E_A| \cdot \Theta(k_A) = \Theta(t_A)$.
By Corollary~\ref{cor:unit+}, $S\cap A$ determines at most $O(t_A
k_A^{\delta-1})= |E_A|\cdot O(k_A^{\delta})$ triangles of any given
area, with (say) $\delta=1/3$. We are now ready to assign triangles
to segments: Assign every triangle $\Delta abc$ determined by $S\cap
A$ to a segment $vw\in E(abc)$ such that (i) every
segment $vw\in E_A$ is assigned to $O(k_A)$ triangles and (ii) every
segment $vw\in E_A$ is assigned to $O(k_A^\delta)$ triangles
of any given area.

\paragraph{Multiplicities.}
Let $T_{vw}(A)$ denote the set of triangles $\Delta abc\subset A$
for which $B(abc)$ is nondegenerate and $\Delta abc$ is assigned to
segment $vw$. We say that the {\em multiplicity} $m(abc)$ of a
triangle $\Delta abc\in T_{vw}(A)$ is the number of triangles in
$T_{vw}(A)$ with the same area as that of $\Delta abc$. For every
$j\in \NN$, let $T_{vw}(A,j)$ denote the set of triangles $\Delta abc\in
T_{vw}(A)$ for which $2^{j-1}\leq m(abc) < 2^j$. Since $vw$ is
assigned to at most  $O(k^\delta_{vw}(A))$ triangles of any given area,
the multiplicity of no triangle can exceed $C k^\delta_{vw}(A)$, for
a sufficiently large constant $C>0$; that is,
$T_{vw}(A,j)=\emptyset$ for $2^j>C k^\delta_{vw}(A)$.

Aggregate the triangles assigned to the same segment $vw$ in various
planes $A\in \A_{vw}$ by letting
$$T_{vw}=\bigcup_{A\in \A_{vw}} T_{vw}(A), \hspace{6mm}\mbox{ and }\hspace{6mm}
T_{vw}(j)=\bigcup_{A\in \A_{vw}} T_{vw}(A,j).$$

\paragraph{Projection to a plane.}
Fix a line segment $vw$. Project $S$ along lines parallel to $vw$ onto a
plane $\pi$ orthogonal to $vw$. The image of $S$ under the projection
is a multiset $\tilde{S}\subset \pi$, where the multiplicity of each
point $\tilde{p}\in \tilde{S}$ is the number of points of $S$ on a
line through $\tilde{p}$ parallel to $vw$. Similarly, the projection of
a subset $Q\subset S$ is a multiset $\tilde{Q}\subset \pi$.
The projection of the line through $vw$ is a single point
$\tilde{v} = \tilde{w}$, and every plane $A$ parallel to $vw$ is
projected into a line $\tilde{A}$ in $\pi$.
An illustration is provided in Fig.~\ref{fig:projection}.

\begin{figure}[htp]
  \begin{center}
\epsfig{file=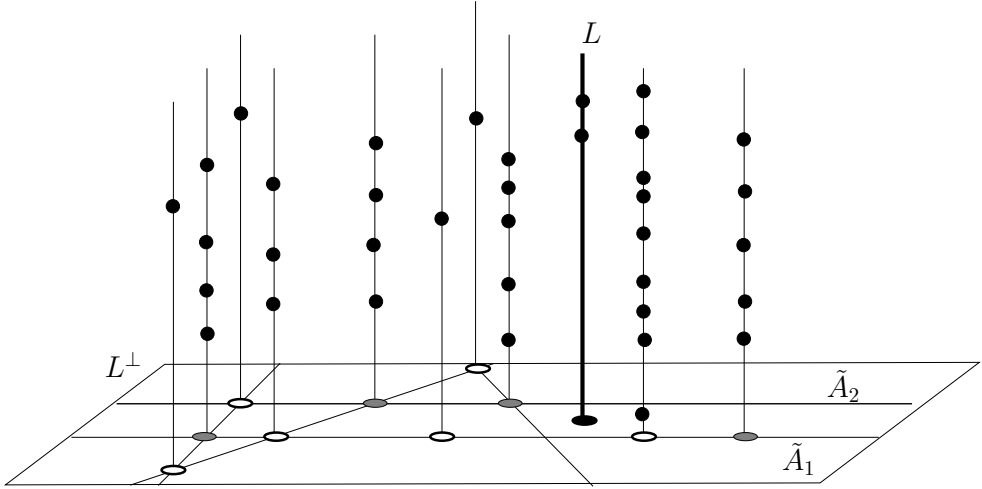,width=13cm,clip=}
 \caption{\small A set $S$ of $n=36$ points in $\RR^3$. Their
 orthogonal projection $\tilde{S}\subset L^\perp$ has 11 points.  The
 multiplicity of every gray point in $L^\perp$ is at least $4$ but
 less than 8, they together form the set $\tilde{S}_{(3)}$. Each of
 the lines $\tilde{A}_1$ and $\tilde{A}_2$ contains two points of
 $\tilde{S}_{(3)}$.
\label{fig:projection}}
\end{center}
 \end{figure}

Sort the points of  $S$ into $\log n$ buckets according to the
multiplicities of their projections: For $r=1,\ldots ,\log n$,
let $S_{(r)}$ denote the set of points $p\in S$ whose projection has
multiplicity at least $2^{r-1}$.
% ; and let $R_{(r)}$ be  the set of points $p\in S$ whose
% projection has multiplicity in the interval $[2^{r-1},2^r)$.
For every plane $A$ parallel to $vw$, we have
$$|S\cap A|=
% \sum_{r=1}^{\log n} |R_{(r)}\cap A| =
% \sum_{r=1}^{\log n} \Theta \left(|\tilde{R}_{(r)}\cap
% \tilde{A}|\cdot 2^r \right) =
\sum_{r=1}^{\log n} O \left( |\tilde{S}_{(r)}\cap \tilde{A}|\cdot 2^r\right). $$

\paragraph{Partition of planes $B(abc)$ according to four parameters.}
We partition the nondegenerate planes parallel to $vw$
into $\log^4 n$ equivalence classes. We will bound the number of
unit volume tetrahedra assigned to $vw$ corresponding to each
equivalence class separately.
For four parameters, $i\in \{0,1,\ldots , \log n\}$, $j\in
\{0,1,\ldots ,\log n \}$, $q\in \{0,1,\ldots , \log n\}$, and $r\in
\{0,1,\ldots , q\}$, we define the set $\B_{vw}(i,j,q,r)$ of
nondegenerate planes $B$ parallel to $vw$ such that there is a unit
volume tetrahedra $abcd$ assigned to $vw$ with $B(abc)=B$ and we
have $2^{i-1}\leq |S\cap A(abc)|<2^i$, $\Delta abc\in T_{vw}(j)$,
and $2^{q-1}\leq 2^r |\tilde{S}_{(r)}\cap B(abc)|<2^q$.

We first give an upper bound on the number of planes in
$\B_{vw}(i,j,q,r)$. There are at most $O(n2^{-i})$ distinct planes
$A\in \A_{vw}$ with $k_{vw}(A)> 2^{i-1}$.
We have assigned at most $O(S\cap A)=O(2^i)$ triangles lying in plane $A$ to
$vw$ by property (i) of the assignment. So there are at most $O(2^{i-j})$ different
triangle areas which are each represented by at most $2^j$ triangles
in $A$ assigned to $vw$.
Triangles $\Delta abc\subset A$ of distinct areas correspond to
distinct planes $B(abc)$ parallel to $A=A(abc)$.  Hence, each such
plane $A$ determines at most $O(2^{i-j})$ parallel
planes $B(abc)$ in $\B_{vw}(i,j,q,r)$, and so $|\B_{vw}(i,j,q,r)|=O(n 2^{-j})$.

We deduce several relation among the parameters $i$, $j$, $q$, and $r$,
assuming $\B_{vw}(i,j,q,r)\neq \emptyset$.
Recall that we assigned $O(2^{\delta i})$ triangles of the same area
to $vw$ in every plane $A\in \A$ with $vw\subset A$ and $k_A\leq
2^i$, by property (ii) of the assignment, and so
$$j \leq \delta i+O(1).$$
Next, we deduce an upper bound on $i$ in terms of $q$ and $r$.
Consider a plane $B\in \B_{vw}(i,j,q,r)$ corresponding to a triangle
$\Delta abc\in T_{vw}(A,j)$ in a plane $A\in \A_{vw}$. By the choice
of $q$ and $r$, at least $2^{q-1}$ points of $S\cap B$ lie on $2^{r-1}$-rich
lines parallel to $vw$. Since we assigned $\Delta abc$ to segment
$vw$ only if $vw$ is not parallel to the $c_0k_{A(abc)}/2=O(2^i)$ richest lines
in $B$, and $B$ contains several $2^{r-1}$-rich lines parallel to $vw$,
we know that an $2^{r-1}$-rich line is not among the $O(2^i)$
richest lines of $B$. By Theorem~\ref{thm:st-83b}, there are at most
$O(2^{2q-3r}+2^{q-r})$ distinct $2^{r-1}$-rich lines in the
plane $B$ that contains at most $2^q$ points. Hence,
the number of $2^{r-1}$-rich lines in any
$B\in \B_{vw}(i,j,q,r)$ is more than $c_02^i/2$. We have
$\Omega(2^i) \leq O(2^{2q-3r}+2^{q-r})$, that is,
$$i\leq \mu, \hspace{6mm} \mbox{where}\hspace{6mm}\mu=\left\{
\begin{array}{lc}
2q-3r +O(1), & \mbox{\rm if } 0\leq r\leq q/2 +O(1),\\
q-r +O(1), & \mbox{\rm if } q/2 -O(1) \leq  r \leq q.
\end{array}\right.$$

\paragraph{Counting incidences in the plane.}
Next, we estimate the number of incidences of the projection points
$\tilde{S}_{(r)}$ and the lines of $\tilde{B}_{vw}(i,j,q,r)$. Since
$|\tilde{\B}_{vw}(i,j,q,r)|=O(n2^{-j})$ and each line spanned by
$\tilde{\B}_{vw}(i,j,q,r)$ contains at least two points of
$\tilde{S}_{(r)}$, Corollary~\ref{cor:st-83} gives
$$ I(\tilde{S}_{(r)},\tilde{\B}_{vw}(i,j,q,r))= O((n2^{-r}\cdot
n2^{-j})^{2/3}+n2^{-r})= O(n^{4/3}2^{-2(j+r)/3}+n2^{-r}). $$
This bound is too weak if $|\tilde{\B}_{vw}(i,j,q,r)|\ll n2^{-j}$, and so
we give another bound for
$I(\tilde{S}_{(r)},\tilde{\B}_{vw}(i,j,q,r))$. Since every
$\tilde{B}\in \tilde{\B}_{vw}(i,j,q,r)$ is an $(2^{q-1}/2^r)$-rich
line in the planar point set $\tilde{S}_{(r)}$, we can apply
Corollary~\ref{cor:st-83b} and we obtain
$I(\tilde{S}_{(r)},\tilde{\B}_{vw}(i,j,q,r))=
O(|\tilde{S}_{(r)}|^2/(2^{q-r})^2+|\tilde{S}_{(r)}|)= O(n^2
2^{-2q}+n2^{-r})$. The combination of these bounds gives
that
\begin{eqnarray*}
I(\tilde{S}_{(r)},\tilde{\B}_{vw}(i,j,q,r))&=& O\left( \min\left(
n^2 2^{-2q}+n2^{-r}, n^{\frac{4}{3}}2^{-\frac{2(j+r)}{3}}+n2^{-r}\right)\right)\\
&=&O\left( \min\left(n^2 2^{-2q},
n^{\frac{4}{3}}2^{-\frac{2(j+r)}{3}}\right)+n2^{-r}\right).
\end{eqnarray*}
In every $B\in \B(i,j,q,r)$, every point of $S\cap B$ corresponds to
$O(2^j)$ unit-volume tetrahedra assigned to $vw$. The number of
points in $B$ is $|S\cap B| = O(\sum_{r=1}^{\log n}
|\tilde{S}_{(r)}\cap B|\cdot 2^r)$. Using the relations we have
observed among the parameters $i$, $j$, $q$, and $r$,  we can sum
these quantities for all $i$, $j$, $q$, and $r$, and get
\begin{eqnarray}
|\U_1(vw)| &=& \sum_{q=1}^{\log n}\sum_{r=1}^q \sum_{i=1}^{\mu}
\sum_{j=1}^{\delta i} \sum_{B\in
\B(i,j,q,r)}|\tilde{S}_{(r)}\cap B|\cdot 2^{j+r}\nonumber\\
&=& O\left( \sum_{q=1}^{\log n}\sum_{r=1}^{q}
\sum_{i=1}^{\mu}\sum_{j=1}^{\delta i}
I(\tilde{S}_{(r)},\tilde{\B}_{vw}(i,j,q,r)) 2^{j+r}\right),\nonumber\\
&=& O\left( \sum_{q=1}^{\log n}\sum_{r=1}^{q}
\sum_{i=1}^{\mu}\sum_{j=1}^{\delta i} \left(\min \left(
n^2 2^{-2q}, n^{\frac{4}{3}} 2^{-\frac{2j+2r}{3}}\right) + n2^{-r}\right) 2^{j+r}\right),\nonumber\\
&=& O\left( \sum_{q=1}^{\log n}\sum_{r=1}^{q}
\sum_{i=1}^{\mu}\sum_{j=1}^{\delta i} \min \left( n^2 2^{j+r-2q},
n^{\frac{4}{3}} 2^{\frac{j+r}{3}}\right) +
n2^j \right),\nonumber\\
&=& O\left( \sum_{q=1}^{\log n}\sum_{r=1}^{q} \sum_{i=1}^{\mu} \min
\left( n^2 2^{\delta i+r-2q}, n^{\frac{4}{3}} 2^{\frac{\delta
i+r}{3}}\right) + n2^{\delta i}
\right),\nonumber\\
&=& O\left( \sum_{q=1}^{\log n}\sum_{r=1}^{q/2} \min \left( n^2
2^{\delta (2q-3r)+r-2q}, n^{\frac{4}{3}} 2^{\frac{\delta
(2q-3r)+r}{3}}\right) + n^{1+\delta}
\right)+\nonumber\\
&& + O\left( \sum_{q=1}^{\log n}\sum_{r=q/2}^{q} \min \left( n^2
2^{\delta (q-r)+r-2q},
n^{\frac{4}{3}} 2^{\frac{\delta (q-r)+r}{3}}\right) + n^{1+\delta} \right),\nonumber\\
&=& O\left( \sum_{q=1}^{\log n}\sum_{r=1}^{q/2} \min \left( n^2
2^{(2\delta -2)q + (1-3\delta)r}, n^{\frac{4}{3}} 2^{\frac{2\delta
q+(1-3\delta)r}{3}}\right) +
n^{1+\delta} \right)+\nonumber\\
&& + O\left( \sum_{q=1}^{\log n}\sum_{r=q/2}^{q} \min \left( n^2
2^{(\delta -2)q +(1-\delta)r},
n^{\frac{4}{3}} 2^{\frac{\delta q +(1-\delta)r}{3}}\right) + n^{1+\delta}\right),\nonumber\\
&=& O\left( \sum_{q=1}^{\log n} \min \left( n^2 2^{\frac{(\delta
-3)q}{2}}, n^{\frac{4}{3}} 2^{\frac{(1+\delta)q}{6}}\right) + \min
\left( n^2 2^{-q},
n^{\frac{4}{3}} 2^{\frac{q}{3}}\right) + n^{1+\delta}\log n\right),\nonumber\\
&=& O\left( n^{\frac{7-\delta}{5-\delta}} + n^{\frac{3}{2}} +
n^{1+\delta}\log^2 n\right) =O\left( n^{3/2} \right).\nonumber
\end{eqnarray}
Summing over all segments $vw$, we obtain
$$|\U_1|=\sum_{(v,w)\in S^2} |\U_1(vw)| =
O(n^2)\cdot O(n^{3/2}) = O(n^{7/2}).$$

\paragraph{Consider $\U_2$.}
For every line $L$ spanned by $S$, consider all planes $A(abc)$
spanned by some triangle $\Delta{abc}$ for which the plane $B(abc)$ is
degenerate and its principal line is $L$.
Project $S$ onto a plane $L^\perp$ orthogonal to $L$. Decompose the
point set $S$ into $\log n$ buckets as described above: For
$r=1,2,\ldots , \log n$, the bucket $\tilde{S}_{(r)}$ contains all
points in $L^\perp$ whose multiplicity is at least $2^{r-1}$.
For a given plane $A$,
observe that the number of $2^{r-1}$-rich lines in $A$ that are
parallel to $L$ is $|\tilde{S}_{(r)}\cap A|$. We know that $|S\cap
A|= O\left(\sum_{r=1}^{\log n} 2^r\cdot |A\cap
\tilde{S}_{(r)}|\right)$. It follows by Jensen's inequality that
$$|S\cap A|^{2+\delta}= O\left( \sum_{r=1}^{\log n} \left(2^{r}\cdot
    |A\cap \tilde{S}_{(r)}|\right)^{2+\delta} \log^{1+\delta} n
\right).$$

For two indices $q\in \{1,2,\ldots , \log n\}$ and $r\in \{1,2,\ldots , q\}$,
let $\A(q,r;L)$ denote the set of planes $A\in \A$
such that $A$ is parallel to $L$, it is incident to at least two
points of $\tilde{A}_{(r)}$, and $2^{q-1}\leq 2^r |A\cap
\tilde{S}_{(r)}|< 2^q$.
The projection $\tilde{A}$ of a plane $A\in \A(q,r;L)$ is a line
incident to at least $(2^{q-1})/2^r = \Omega(2^{q-r})$ points in
bucket $\tilde{S}_{(r)}$. We clearly have $|\tilde{S}_{(r)}| = O(n2^{-r})$,
and so the number of such lines can be bounded by
Theorem~\ref{thm:st-83b}:
$$|\A(q,r;L)| = O \left( \frac{(n2^{-r})^2}{(2^{q-r})^3} +
\frac{n2^{-r}}{2^{q-r}}\right) = O \left( n^2 2^{r-3q} + n2^{-q}
\right).$$
For a fixed line $L$, denote by $\A(L)$ the set of planes
spanned by $S$ and containing $L$.
By the fact that the number of unit area triangles in a plane $A
\in  \A(L)$ is $O(|S\cap A|^{2+\delta})$,
the number of triangles $\Delta abc$ whose degenerate plane $B(abc)$
contains $L$ is
\begin{eqnarray}
|T_2(L)|&= &\sum_{q=1}^{\log n}\sum_{r=1}^{q} |\A(q,r;L)|\cdot
O(2^{(2+\delta)q}\log^{1+\delta}n) =\sum_{q=1}^{\log
n}\sum_{r=1}^{q}
O \left( n^2 2^{r+ (\delta-1)q} + n 2^{(1+\delta)q}
\right)\log^{1+\delta} n \nonumber\\
&=&\sum_{q=1}^{\log n} O \left( n^2 2^{\delta q} + qn
2^{(1+\delta)q}\right)\log^{1+\delta} n
=O \left(n^{2+\delta} \log^{2+\delta} n
\right) \nonumber .
\end{eqnarray}

By Theorem~\ref{thm:st-83b},  for every $M=O(\sqrt{n})$
there are at most $O(n^2/M^3)$ distinct $M$-rich lines. Moreover,
the total number of incidences on these lines is at most
$O(n^2/M^2)$. Therefore, the number of unit volume tetrahedra $abcd
\in \U_2$ whose plane $B(abc)$ is degenerate and its principal line
$L$ is $M$-rich is bounded by $O(|T_2(L)| \cdot n^2/M^2)=
O\left(n^{4+\delta} \log^{2+\delta} n/M^2\right)$.
On the other hand, the number of unit tetrahedra $abcd \in \U_2$ whose
plane $B(abc)$ is degenerate and its principal line $L$ contains at most
$M$ points, and so $|S\cap B(abc)|=O(M)$, is at most $O(n^3 M)$.
We balance the two upper bounds by choosing $M:=n^{(1+\delta)/3}\log^{(2+\delta)/3} n$,
and then we have $|\U_2| = O(n^{(10+\delta)/3}\log^{(2+\delta)/3} n)$.

Recalling that $\delta\leq 1/3$, we have proved that
$$|\U| = |\U_1|+ |\U_2| =
O \left(n^{7/2} \right)+O\left(n^{(10+\delta)/3}\log^{(2+\delta)/3}
n\right)= O\left(n^{7/2}\right). $$
\end{proof}

\section{The number of distinct volumes\label{S-distinct}}

In 1982, Erd\H os, Purdy, and Straus \cite{eps-82} considered the
analogue of the problem of distinct triangle areas in the plane to
higher dimensions and posed the following problem:
Let $S$ be a set of $n$ points in $\RR^d$ not all in one hyperplane.
What is the minimal number $g_d(n)$ of distinct volumes of
full-dimensional simplices with vertices in $S$?

By taking $d$ sets of about $n/d$ equally spaced points on parallel
lines through the vertices of a $(d-1)$-simplex, one gets $g_d(n)
\leq \lfloor \frac{n-1}{d} \rfloor$. Erd\H os, Purdy, and Straus
conjectured that equality holds at least for sufficiently large $n$
(see also \cite{cfg-91}).

In this section, we  give a first linear lower bound on the number of
full-dimensional simplices in $\RR^d$ (including tetrahedra in
3-space) determined by $n$ points that span  $\RR^d$.
Some remarks are in order:
We have mentioned earlier Erd\H{o}s's distinct distance problem that asks for
the minimum number $t(n)$ of distinct inter-point distances among $n$
points in the plane. The currently known best lower bound for $t(n)$,
due to Katz and Tardos~\cite{kt-04}, in fact gives a lower bound on
the maximum number of interpoint distances measured {\em from a single
point}.  For triangle areas in the plane, we have a similar
phenomenon. Pinchasi proved that every set $S$ of $n$ noncollinear
points in the plane contains a point pair $a,b\in S$ such that the
points of $S$ determine $\lfloor (n-1)/2\rfloor$ distinct (nonzero)
distances measured from the line $ab$, and so there are at least $\lfloor
(n-1)/2\rfloor$ triangles with distinct (positive) areas \cite{p-07}.

Our proof follows a similar path. We show that $n$ points in
$d$-space, not all on a hyperplane, determine $\Omega(n)$
full-dimensional simplices  with distinct volumes that share a common
face. We start with a useful lemma about triangle areas in the plane:

\begin{lemma}\label{lem:origin}
  Let $S$ be a set of $n$ points in the plane and let $p_1 \in S$ be a
  point such that each line $p_1q$, $q \in S\setminus \{p_1\}$, is
  incident to no other points in $S$. Then there is a point $p_2\in
  S\setminus \{p_1\}$ such that the triangles $\Delta p_1 p_2 q$ with
  $q \in S\setminus \{p_1,p_2\}$ determine $\Omega(n)$ distinct areas.
\end{lemma}
\begin{proof}
  Let $t\in \NN$ be the smallest integer such that for every $p_2\in
  S\setminus \{p_1\}$, the triangles $S$ with a common side
  $p_1p_2$ have at most $t$ distinct areas. For every
  fixed $p_2\in S\setminus\{p_1\}$, al points $q\in S\setminus \{p_1,p_2\}$ where
  the area of the triangle $\Delta p_1p_2q$ is one of the at most $t$ possible areas
  must lie on one of at most $2t$ lines parallel to $p_1p_2$. Hence for every
  $p_2\in S \setminus \{p_1\}$, all the remaining $n-2$ points of $S$
  must lie on one of at most $2t$ lines.  Since the directions of the
  lines $p_1p_2$, $p_2\in S\setminus \{p_1\}$, are all distinct, we
  obtain a total of $\ell\leq (2t+1)(n-1)$ distinct lines (the lines
  $p_1p_2$ are also counted here, accounting for collinear triples $(p_1,p_2,q)$).
  Consequently there are $n(n-1)$ point-line incidences between $S$ and
  these $\ell$ distinct lines.

  By the Szemer\'edi-Trotter Theorem, the number of incidences is at
  most $O(n^{2/3}\ell^{2/3}+n+\ell)=$ $O(n^{4/3}t^{2/3}+n+nt)$. This
  implies that $t=\Omega(n)$.
\end{proof}

\begin{theorem}\label{thm:distinct}
  The tetrahedra determined by a set $S$ of $n$ points in $\RR^3$, not
  all on a plane, have at least $\Omega(n)$ distinct volumes.
  Furthermore, $S$ spans a triangle $F$ such that $S$ determines
  $\Omega(n)$ tetrahedra of distinct volumes, having a common face
  $F$.
\end{theorem}
\begin{proof}
  Let $S$ be a set of $n$ points in $\RR^3$, not all on a plane.
  Consider the two cases from Lemma~\ref{lem:b}. First assume that
  there are $\beta_3 n$ points in a plane $A$. Since $S$ is not
  contained in a plane, we may assume
  further that $S$ spans the plane $A$ (not all points of $S\cap A$
  are collinear).
  By the planar result mentioned above, the points $S\cap A$ span at least
  $\Omega(|S\cap A|)= \Omega(n)$ triangles with distinct areas and
  with a common side
  $pq$. Let $r \in S\setminus A$. The triangles in $A$ together with
  $r$ determine $\Omega(n)$ tetrahedra of distinct volumes with a
  common face $\Delta pqr$.

  Next, assume that the triples of $S$ span at least $\gamma_3 n^3$
  distinct planes. There is a segment $p_0p_1$ that participates
  in at least $\gamma_3 n^3/ {n\choose 2} =\Omega(n)$ such triples.
  Let $L$ be the line through $p_0p_1$; and let $\A$ be the set of
  planes determined by $S$ and containing $p_0p_1$. Let $S_0$,
  $S_0\subset S\setminus L$, be a set of points that contains exactly one
  point from each $A\setminus L$, $A\in \A$.

  Project the point set $S_0$ onto a plane $\pi$ orthogonal to $L$.
  Refer to Fig.~\ref{fig:distinct}.  Let $\tilde{S}_0\subset \pi$
  denote the projection points (and the projection of $L$ onto $\pi$
  is a point $\tilde{p}_1$).  By construction, $\tilde{S}_0\subset
  \pi$ contains $\Omega(n)$ points such that each line
  $\tilde{p}_1\tilde{q}$, $\tilde{q}\in \tilde{S}_0$, is incident to a
  unique point of $\tilde{S}_0$.
\begin{figure}[htbp]
\centerline{\epsfxsize=4.4in \epsffile{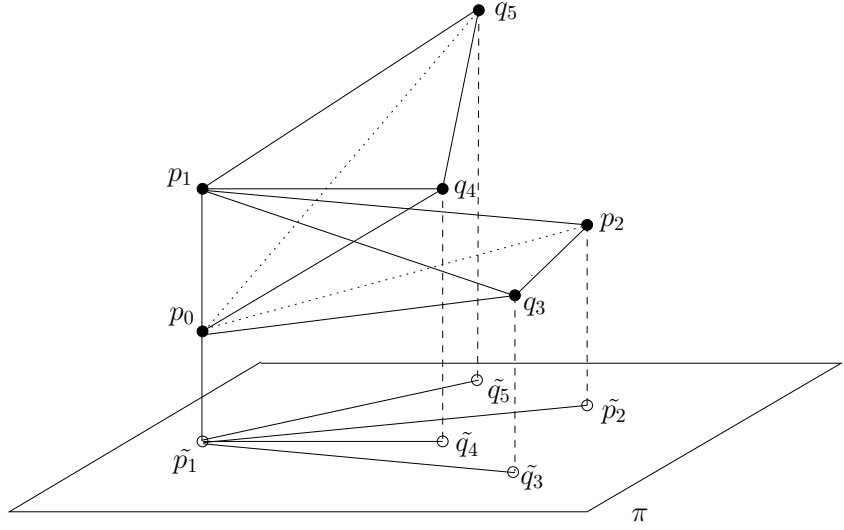}}
\caption{\small Tetrahedra $T(p_0p_1p_2q_3)$, $T(p_0p_1p_2q_4)$ and
$T(p_0p_1p_2q_5)$ have distinct volumes and a common face $\Delta
p_0p_1p_2$.}
\label{fig:distinct}
\end{figure}
  By Lemma~\ref{lem:origin}, there is a point $\tilde{p}_2\in
  \tilde{S}_0$ such that the triangles $\Delta
  \tilde{p}_1\tilde{p}_2\tilde{q}$ with $\tilde{q}\in
  \tilde{S}_0\setminus \{\tilde{p}_1,\tilde{p}_2\}$ determine
  $\Omega(n)$ distinct areas in the plane $\pi$.  For every triangle
  $\Delta \tilde{p}_1\tilde{p}_2\tilde{q}$,
  the volume of the tetrahedron $T(p_0 p_1 p_2 q)$ can be expressed as
  (here we use again the fact that if the vertices of a tetrahedron
  are on three parallel lines, one can shift the single points along
  the corresponding lines while the volume remains the same):
  $${\rm vol}(T(p_0 p_1 p_2 q)) = \frac{1}{3}{\rm area} (\Delta
  \tilde{p}_1\tilde{p}_2\tilde{q}) \cdot |p_0p_1| .$$
  Hence $S$
  contains $\Omega(n)$ tetrahedra of distinct volumes sharing a common
  face $\Delta p_0p_1p_2$.
\end{proof}

Theorem~\ref{thm:distinct} readily generalizes for full-dimensional
simplices in any dimension $d\in \NN$.

\begin{theorem}\label{thm:ddistinct}
  For every $d\in \NN$, the full-dimensional simplices determined by a
  set $S$ of $n$ points in $\RR^d$, not all on a hyperplane, have at
  least $\Omega(n)$ distinct volumes.  Furthermore, $S$ spans a
  $(d-1)$-dimensional simplex $F$ such that $S$ determines $\Omega(n)$
  full-dimensional simplices of distinct volumes, having a common face
  $F$.
\end{theorem}
\begin{proof}
  We proceed by induction on $d\in \NN$. The theorem is obvious for
  $d=1$ and it was confirmed by Burton and Purdy for
  $d=2$~\cite{bp-79}.  Let $d\geq 3$, and assume that for every
  $d'<d$, there is a constant $c(d')>0$ such that any $m$-element
  point set in $\RR^{d'}$, not all in a $(d'-1)$-dimensional affine
  subspace, determines at least $c(d')m$ \ $d'$-dimensional
  simplices of distinct volumes having a common $(d'-1)$-dimensional
  face.

  Let $S$ be a set of $n$ points in $\RR^d$, not all on a hyperplane.
  Consider the two cases from Lemma~\ref{lem:b}: First assume that
  there are $\beta_d n$ points in a hyperplane. These points
  span a $d'$-dimensional subspace $A$, for some $1\leq d'<d$.  By the
  induction hypothesis, the points $S\cap A$ determine $c(d')\beta_d
  n=\Omega(n)$ \ $d'$-dimensional simplices of distinct volumes having
  a common $(d'-1)$-dimensional face. Since $S$ is not contained in a
  hyperplane, there is a set $T_A\subset S\setminus A$ of $d-d'$
  points such that $T_A$ and $A$ together span the entire space $\RR^d$.
  Adding the vertices $T_A$ to each of the $\Omega(n)$ \ $d'$-dimensional
  simplices of distinct volumes in $A$, we obtain $\Omega(n)$
  full-dimensional simplices of distinct volumes having a
  common face in $\RR^d$.

  Next, assume that the $d$-tuples of $S$ span at least $\gamma_d n^d$
  distinct hyperplanes. There is a $(d-1)$-tuple $F \subset S$ that
  participates in at least $\gamma_d n^d/ {n\choose d-1} =\Omega(n)$
  such $d$-tuples spanning distinct hyperplanes.  Let $A$ be the
$(d-2)$-dimensional affine subspace
  spanned by $F$; and let $\L$ be the set of $(d-1)$-dimensional
  affine subspaces determined by $S$ and containing $F$. Let $S_0$,
  $S_0\subset S\setminus A$, be a set that contains exactly one point from
  $B\setminus A$, for each $B\in \L$.

  Project the point set $S_0$ onto a plane $\pi$ orthogonal to $A$.
  Refer to Fig.~\ref{fig:distinct}.  Let $\tilde{S}_0\subset \pi$
  denote the projection points (and the projection of $A$ onto $\pi$
  is a point $\tilde{p}_1$).  By construction, $\tilde{S}_0\subset
  \pi$ contains $\Omega(n)$ points such that each line
  $\tilde{p}_1\tilde{q}$, $\tilde{q}\in \tilde{S}_0$, is incident to a
  unique point of $\tilde{S}_0$.

  By Lemma~\ref{lem:origin}, there is a point $\tilde{p}_2\in
  \tilde{S}_0$ such that the triangles $\Delta
  \tilde{p}_1\tilde{p}_2\tilde{q}$ with $\tilde{q}\in
  \tilde{S}_0\setminus \{\tilde{p}_1,\tilde{p}_2\}$ determine
  $\Omega(n)$ distinct areas in the plane $\pi$. For every triangle
  $\Delta \tilde{p}_1\tilde{p}_2\tilde{q}$, the
  volume of the full-dimensional simplex $T(F\cup \{p_2,q\})$ spanned by
  $F\cup \{p_2,q\}$ can be expressed as:
  $${\rm vol}(T(F \cup \{p_2,q\})) = \frac{2!(d-2)!}{d!}\cdot {\rm
    vol}(\Delta \tilde{p}_1\tilde{p}_2\tilde{q}) \cdot {\rm vol}(F).$$
  Hence $S$ contains $\Omega(n)$ full-dimensional simplices of
  distinct volumes sharing a common face $F \cup p_2$.
\end{proof}

Our proof for Theorem~\ref{thm:ddistinct} crucially relies on the fact that we
build {\em full-dimensional} simplices of distinct volumes (i.e., $k=d$).
No tight bound is known for the order of magnitude of the minimum number of
$k$-dimensional simplices of distinct volumes determined by $n$ points
in $\RR^d$, not all in a hyperplane, for any $1\leq k<d$. In
particular, for $k=1$ this is Erd\H{o}s's celebrated distinct distance
problem in $\RR^d$ we mentioned earlier.

\end{document}